\documentclass{amsart}
\numberwithin{equation}{section}
\usepackage{latexsym}
\usepackage{amsmath}
\usepackage{amssymb}
\usepackage{mathrsfs}
\usepackage{graphicx}
\usepackage{ifthen}
\usepackage{comment}
\usepackage[T1]{fontenc}
\usepackage[latin1]{inputenc}

\numberwithin{equation}{section}

\newtheorem{definition}{Definition}[section]
\newtheorem{theorem}[definition]{Theorem}
\newtheorem{lemma}[definition]{Lemma}
\newtheorem{corollary}[definition]{Corollary}
\newtheorem{proposition}[definition]{Proposition}
\newtheorem{remark}[definition]{Remark}

\newtheorem{example}[definition]{Example}

\usepackage{latexsym}
\usepackage{amsmath}
\usepackage{amssymb}

\newcommand{\R}{{\mathbb R}}

\newcommand{\ve}{{\varepsilon}}

\def\ve{\varepsilon}

\def\d{\partial}

\def\F{\mathcal{F}}

\def\f{\varphi}

\def\R{\mathbb{R}}
\def\E{\mathbb{E}}

\begin{document}

\title{Invariant Measure for Stochastic Functional Differential Equations in Hilbert Spaces}

\author{Oleksandr Misiats}
\address{Department of Mathematics\\
Virginia Commonwealth University, Richmond, VA, USA}
\email{omisiats@vcu.edu}

\author{Viktoriia Mogylova}
\address{Department of Physics and Mathematics\\
Igor Sikorsky Kyiv Polytechnic Institute, Ukraine}
\email{mogylova.viktoria@gmail.com}

\author{Oleksandr Stanzhytskyi}
\address{Department of Mathematics, Taras Shevchenko National University of Kyiv, Ukraine}
\email{ostanzh@gmail.com}


\begin{abstract}
In this work we study the long time behavior of nonlinear stochastic functional-differential equations in Hilbert spaces. In particular, we start with establishing the existence and uniqueness of mild solutions. We proceed with deriving a priory uniform in time bounds for the solutions in the appropriate Hilbert spaces. These bounds enable us to establish the existence of invariant measure based on Krylov-Bogoliubov theorem on the tightness of the family of measures.  Finally, under certain assumptions on nonlinearities, we establish the uniqueness of invariant measures. 
\end{abstract}

\maketitle

\vspace{-.2cm}
\section{Introduction}
In this work we study the asymptotic behaviour of the solutions of stochastic functional-differential equations. In a bounded domain, the equation reads as
\begin{align}\label{MainB}
& du = [Au + f(u_t)] dt + \sigma(u_t) d W(t) \text{ in } D, t>0; &\\
\nonumber & u(t,x)= \phi(t,x), t \in [-h,0), u(0,x) = \f_0(x)   \text{ in } D;& \\
\nonumber & u(t,x) = 0, x \in \partial D, t \geq 0.&
\end{align}
The corresponding problem in the entire space has the form
\begin{align}\label{MainUB}
& du = [Au + f(u_t)] dt + \sigma(u_t) d W(t) \text{ in } \R^d, t>0; &\\
\nonumber & u(t,x)= \phi(t,x), t \in [-h,0), u(0,x) = \f_0(x)   \text{ in } \R^d.&
\end{align}
Here $A$ is an elliptic operator
\begin{equation}\label{operator}
A = A(x) = \sum_{i,j=1}^{d} a_{ij}(x) \frac{\d^2}{\d x_i \d x_j} + \sum_{i=1}^{d} b_i(x) \frac{\d}{\d x_i} + c(x),
\end{equation}
the interval $[-h,0]$ is the interval of delay, and $u_t = u(t+ \theta)$ with $\theta \in [-h,0)$.

Functional differential equations of types (\ref{MainB}) and (\ref{MainUB}) are mathematical models of processes, the evolution of which depends on the previous states. Such models are widely used in population dynamics, electrical engineering \cite{Jason}, chemical engineering \cite{WanFan} et. al. The classic results for deterministic functional-differential equations in finite dimensional spaces can be found in \cite{Hal77} and references therein. Stochastic functional differential equation in finite dimensions have be studies extensively as well. In particular, the existence of invariant measures for stochastic ordinary differential equations was established in \cite{ButSch} and \cite{HaiMatSch11}. The work \cite{IvaKazSwi} addressed the stochastic stability, as well as various applications of stochastic delay equations in finite dimensions. 

The results on functional differential equations in infinite dimensions are significantly more sparse. One example of analysis and applications of functional partial differential equations may be found in \cite{OmaGur}. In this work, the authors study the nonlocal reaction-diffusion model of population dynamics. They establish the existence of time stationary solution and show that all other solutions converge to it. 

The results on stochastic functional differential equations include \cite{Tan02} and \cite{CuiYanSun}, which establish the existence of solutions and their stability. Stochastic differential equation of neutral type were studied in \cite{SamMahSta08} and \cite{KenStaTsu}. The work \cite{StaMogTsu19} established the comparison principle for such equations.

The main goal of the present work is to establish the existence and uniqueness of invariant measures for the equations \eqref{MainB} and \eqref{MainUB} based on Krylov-Bogoliubov theorem on the tightness of the family of measures \cite{KryBog}. More precisely, we will use the compactness approach of Da Parto and Zabczyk \cite{DapZab92}, which involves the following key steps:
\begin{itemize}
  \item[{[i]}] Establishing the existence of a Markovian solution of (\ref{MainB}) or (\ref{MainUB}) in a certain functional space, in which the corresponding transition semigroup is Feller;
  \item[{[ii]}] Showing that the semigroup $S(t)$ generated by $A$ is compact;
  \item[{[iii]}] Showing that the corresponding equation with a suitable initial condition has a solution, which is bounded in probability.
\end{itemize}
This approach was used in establishing the existence of invariant measure for a large class of stochastic nonlinear partial differential equations without delay, e.g. \cite{DapZab96}, \cite{Chow}, \cite{AssMan01}, \cite{Cer99}, \cite{MisStaYip2}, \cite{MisStaYip3} and references therein. 
For functional differential equations in finite dimensions, the approach above was used in \cite{Cho-Mic}. In this work, the author established the existence of an invariant measure in $\R^d \times L^2(-h,0;\R^d)$. In contrast, for stochastic partial differential equations, the natural phase space for the mild solutions of (\ref{MainUB}) is
$L^2_\rho(\R^d) \times L^2(-h,0,L^2_{\rho}(\R^d))$, where $L^2_\rho(\R^d)$ is a weighted space. The equations of type (\ref{MainB}) and (\ref{MainUB}) were studied in the space $C(-1,0,L^2_{\rho}(\R^d))$, which is a significantly easier problem \cite{SamMahSta08}, \cite{StaMogTsu19}, \cite{Tan02}.
In these spaces the authors studied the conditions for the existence and uniqueness of the solution, as well as their Markov's and Feller properties. However, in order to apply the compactness approach one needs to work in $L^2_\rho(\R^d) \times L^2(-h,0,L^2_{\rho}(\R^d))$, which is done in this work. We also establish the existence and uniqueness of the stationary solution, and the convergence of other solutions to it in square mean, which is the stochastic analog of the main result of \cite{OmaGur}. 

This paper is structured as follows. In Section 2 we introduce the notation and formulate the main results. Section 3 is devoted to the proof of the existence of invariant measure, as well as an example of application of this result to integral-differential equations. Section 4 establishes the uniqueness of invariant measure, and the convergence to the stationary solution.

\section{Preliminaries and Main Results.} \label{Ch2} Throughout the paper, the domain $D$ is either a bounded domain with $\d D$ satisfying the Lyapunov condition, or $D = \R^d$. Denote
\begin{equation}\label{weight}
\rho(x):= \frac{1}{1+|x|^r}
\end{equation}
where $r>d$ if $D = \R^d$ and $r = 0$ (i.e. no weight) for bounded $D$. We introduce the following spaces:
\begin{align}\label{spaces}
  B_0^{\rho} & := L^2_{\rho}(D); \\
 \nonumber B_1^{\rho} & := L^2(-h,0,L^2_{\rho}(D)); \\
\nonumber B^{\rho} & := B_0^{\rho} \times B_1^{\rho};\\
\nonumber H & := L^2(D).
\end{align}
The coefficients $a_{ij}$ of the operator $A$ defined in (\ref{operator}) are Holder continuous with the exponent $\beta \in (0,1)$, symmetric, bounded and satisfying the elipticity condition
\[
\sum_{i,j=1}^{d} a_{i,j} \eta_i \eta_j \geq C_0 |\eta|, \ \eta \in \R^d
\]
for some $C_0>0$. The coefficients $b_i$ and $c$ are also bounded and Holder continuous with some positive Holder exponent.

If $D$ is bounded, we impose homogeneous Dirichlet boundary conditions on $\partial D$. In this case,
\[
D(A) = H^2(D) \cap H_0^1(D).
\]
If $D = \R^d$, then $D(A) = H^2(\R^d).$
Denote $G(t,x,y)$ to be the fundamental solution (or the Green's function in the case of bounded $D$) for $\frac{\d}{\d t} - A$. It follows from, e.g., \cite{LadSolUra}, p. 468, that there are positive constants $C_1(T), C_2(T) > 0$ such that
\begin{equation}\label{NashAr}
0 \leq G(t,x,y) \leq C_1(T) t^{-d/2} e^{-C_2(T) \frac{|x-y|^2}{t}}
\end{equation}
for $t \in [0,T]$ and $x,y \in D$. Note that in (\ref{NashAr}), $C_1$ and $C_2$ depend not only on $T$, but on the constants $C_0$, $d$, $T$, maximum values of the coefficients of $A$, and the Holder constants. If the operator is in the divergence form $A u = \rm{ div }(a \nabla u)$, the estimates are of a different type, see e.g. [17], namely
\begin{equation}\label{GFest}
g_1(t,x-y) \leq G(t,x,y) \leq g_2(t, x-y)
\end{equation}
where
\[
g_i(t,x) = K(C_0,d) t^{-d/2} e^{-K(C_0,d) \frac{|x|^2}{t}}, \ i=1,2, t \geq 0, \ x, y \in \R^d.
\]
In this case, in contrast with (\ref{NashAr}),  the constant $K(C_0,d)$ is independent on $t$.

\begin{lemma}\label{L2}
For all $T>0$ there exists a positive $C(r,T)>0$ such that
\begin{equation}\label{weight2}
    \int_{D} G(t,x,y) \rho(y) dy \leq C(r,T) \rho(x), t \in [0,T].
\end{equation}
\end{lemma}
\begin{proof}
Note that the weight (\ref{weight}) satisfies 
\begin{equation}\label{weight1}
\frac{\rho(x)}{\rho(y)} \leq C(r)(1+|x-y|^r)
\end{equation}
for some $C(r)>0$. Thus
\[
\int_{D} G(t,x,y) \rho(y) dy \leq C(r) \int_{D}G(t,x,y) \rho^{-1}(x-y) \rho(x) dy \leq
\]
\[
\leq C(r) C_1(T)  \int_{\R^d} t^{-d/2} e^{-C_2(T) \frac{|y|^2}{t}} (1+|y|^r) dy \rho(x) \leq C(r,T) \rho(x).
\]
\end{proof}
Define
\begin{equation}\label{def_semi}
    (S(t) \f)(x):= \int_{D}G(t,x,y) \f(y) \, dy, t>0, x \in D, \f \in L^2(D), \text{ and } S(0) = I,
\end{equation}
which is a semigroup on $L^2(D)$ with generator $A$. Then for all $\f \in L^2(D)$ and for $t \in [0,T]$ by Lemma \ref{L2} we have
\begin{equation}\label{ast}
    \|S(t) \f \|_{\rho}^2 = \int_{D} \left(\int_{D} G(t,x,y) \f(y) dy \right)^2 \rho(x) dx \leq
\end{equation}
\[
\leq \int_{D} \rho(x) \left(\int_{D} G(t,x,y) dy\right)\left(\int_{D} G(t,x,y) \f^2(y) dy\right) \leq
\]
\[
\leq C \int_{D} \left(\int_{D} G(t,x,y) \frac{\rho(x)}{\rho(y)} dx\right) \rho(y) \f^2(y) dy \leq C_{\rho}(T) \|\f\|_{\rho}^2.
\]
The above estimate allows the semigroup $S(t)$ to be extended to a linear map from $B_{0}^{\rho}$ to itself. Since $L^2(D)$ is dense in $B_0^{\rho}$, $S(t)$ is strongly continouos in $B_0^{\rho}$.

Let $\sum_{i=1}^{\infty} a_i < \infty$, and $e_n$ be orthonormal basis in $H$, such that $e_n \in L^{\infty}(D)$ and $\sup_{n}\|e_n\|_{L^{\infty}(D)}<\infty$. Introduce the operator $Q \in \mathcal{L}(H)$ such that $Q$ is non-negative, $Tr(Q)< \infty$, $Q e_n  = a_n e_n$.
Let $(\Omega, \mathcal{F},P)$ be a complete probability space. Introduce
\[
W(t):= \sum_{i=1}^{\infty} \sqrt{a_i}\beta_i(t) e_i(x), \ t\geq 0,
\]
which is a $Q$-Wiener process on $t \geq 0$ with values in $L^2(Q)$. Here $\beta_i(t)$ are standard, one dimensional, independent Wiener processes.  Also let $\{F_{t}, t \geq 0\}$ be a normal filtration satisfying
\begin{itemize}
    \item $W(t)$ is $\mathcal{F}_t$-measurable;
    \item $W(t+h)-W(t)$ is independent of $\mathcal{F}_t$ $\forall h \geq 0, t\geq 0$.
\end{itemize}
Denote $U = Q^{\frac{1}{2}}(H)$. From \cite{ManZau99}, Lemma 2.2, $U \in L^{\infty}(D)$. Following \cite{ManZau99} introduce the multiplication operator $\Phi:U \to B_0^{\rho}$ as follows: for fixed $\f \in B_0^{\rho}$, let $\Phi(\psi): = \f \psi, \ \psi \in U$. Since $\f \in B_0^{\rho}$ and $\f \in L^{\infty}(D)$, the operator is well defined and hence $\Phi \circ Q^{1/2}: L^2(D) \to B_0^{\rho}$ defines a Hilbert-Schmidt operator. The operator $\Phi$ is also a Hilbert-Schmidt operator satisfying
\[
\|\Phi \circ Q^{1/2} \|^2_{\mathcal{L}_2} = \sum_{n=1}^{\infty} \|\Phi \circ Q^{1/2} e_n\|^2_{B_0^\rho} = \sum_{n=1}^{\infty} a_n  \int_{D} \f^2(x) e_n^2(x) \rho(x) dx \leq
\]
\[
\leq Tr(Q)  \sup_{n} \|e_n\|_{\infty}^2 \|\f\|^2_{\rho},
\]
where $Tr Q = \sum_{n=1}^{\infty} a_n = a$. Hence if $\Phi: \Omega \times [0,T] \to \mathcal{L}(U, B_0^{\rho})$ is a predictable process satisfying
\[
\E \int_{0}^{T} \| \Phi \circ Q^{1/2}\|_{\mathcal{L}_2}^2 ds < \infty,
\]
following \cite{DapZab92} we can define
\[
\int_0^t \Psi(s) d W(s) \in B_0^{\rho}
\]
with the following expansion
\[
\int_0^t \Psi(s) d W(s) = \sum_{i=1}^{\infty} \sqrt{a_i} \int_0^t \Phi(s, \cdot) e_i(\cdot) d \beta_i(s).
\]
Furthermore,
\begin{equation}\label{Noise_est}
\E \left\| \int_0^t \Psi(s) d W(s) \right\|_{\rho}^2 \leq a \sup_n \|e_n\|_{\infty}^2  \int_0^t \E \|\Psi(s,\cdot)\|^2_{B_0^{\rho}} ds.
\end{equation}

{\bf Assumptions on nonlinearities.} Assume $f$ and $\sigma$ satisfy the following conditions:
\begin{enumerate}
    \item[{[i]}] The functionals $f$ and $\sigma$ map $B_1^{\rho}$ to $B_0^{\rho}$;
    \item[{[ii]}] There exists a constant $L>0$ such that
    \[
    \|f(\f_1) - f(\f_2)\|_{B_0^{\rho}} + \|\sigma(\f_1) - \sigma(\f_2)\|_{B_0^{\rho}} \leq L \|\f_1 - \f_2\|_{B_1^{\rho}}
    \]
    for any $\f_1, \f_2 \in B_1^{\rho}$.
\end{enumerate}
\begin{definition}
An $\mathcal{F}_t$ measurable random process $u(t,\cdot) \in B_0^{\rho}$ is a mild solution of (\ref{MainB}) or (\ref{MainUB}), if
\begin{equation}\label{MildSOl}
    u(t,\cdot) = S(t) \f(0,\cdot) + \int_0^t S(t-s) f(u_s) ds + \int_0^t S(t-s) \sigma(u_s) d W(s)
\end{equation}
where
\[
u(0,\cdot) = \f(0, \cdot) \in B_0^{\rho}, \ u(t, \cdot) = \f(t, \cdot) \in B_1^{\rho}, t \in [-h,0].
\]
\end{definition}

\begin{theorem}\label{Th:1}(Existence and uniqueness). Suppose $f$ and $\sigma$ satisfy the conditions [i] and [ii], and $\f(t,\cdot)$  is an $\mathcal{F}_0$ measurable random process for $t \in [-h,0]$, which is independent of $W$ and such that
\[
\E \|\f(0,\cdot)\|_{B_0^{\rho}}^p < \infty \text{ and } \E \|\f(\cdot,\cdot)\|_{B_1^{\rho}}^p < \infty, p \geq 2.
\]
Then there exists a unique mild solution of (\ref{MainB}) (or \ref{MainUB}) on $[0,T]$, and
\begin{equation}\label{cont dep}
\E \|y(t)\|_{B^\rho}^p \leq K(T) (1+\E \|y(0)\|_{B^\rho}^p), \ t \in [0,T].
\end{equation}
\end{theorem}

\begin{theorem}\label{Th:2}(Continuous dependence on the initial data) Let $\phi \in B_1^{\rho}$, $\phi(0, \cdot) \in B_0^{\rho}$,$\phi_1 \in B_1^{\rho}$, $\phi_1(0, \cdot) \in B_0^{\rho}$ be two initial sets of data of two solutions
\[
y(t) = y(t,\phi) = \begin{pmatrix}
u(t,\phi) \\
u_t(\phi)
\end{pmatrix}, y_1(t) = y(t,\phi_1) = \begin{pmatrix}
u(t,\phi_1) \\
u_t(\phi_1)
\end{pmatrix}
\]
respectively. Then under the conditions of Theorem \ref{Th:1} there exists a constant $C(T)$ such that
\begin{equation}\label{cont}
    \sup_{t \in [0,T]} \E \|y(t) - y_1(t)\|_{B^{\rho}}^2 \leq C(T) \E \|\phi(t) - \phi_1(t)\|_{B^{\rho}}^2.
\end{equation}
\end{theorem}

The following proposition shows the that the solution $u(t,\cdot)$ has continuous trajectories.
\begin{proposition}\label{Prop:1}
Let $u(t,\cdot)$ be a mild solution of (\ref{MainB}) or (\ref{MainUB}). Then, under the conditions of Theorem \ref{Th:1} $u_t$ is continuous at $t=0$ in probability with respect to the norm $\|\cdot\|_{B_1^{\rho}}$, i.e.
\[
\|u_t - u_0\|_{B_1^\rho}^2 = \int_{-h}^{-t} \E \|u(t+ \theta) - \f(\theta)\|_{B_0^{\rho}}^2 d\theta \to^P 0, t \to 0.
\]
\end{proposition}
\begin{proof}
\[
\E \|u_t - u_0\|^2_{B_1^{\rho}} \leq \int_{-h}^{-t} \E \|\f(t+\theta) - \f(\theta)\|^2_{B_0^{\rho}} d\theta
+ \int_{-t}^{0} \E \|u(t+\theta) - \f(\theta)\|^2_{B_0^{\rho}} d\theta.
\]
The convergence of the first term to $0$ follows from the density of $C([-h,0], B_0^{\rho} \times L^2(\Omega))$ in $L^2([-h,0], B_0^{\rho} \times B_0^{\rho} \times L^2(\Omega))$. The second term converges to zero as $t \to 0 $ since the integrand is bounded.

\end{proof}

Denote $B_b(B^{\rho})$ to be the Banach space of bounded real Borel functions from $B^{\rho}$ to $\R$, and $C_b (B^{\rho})$ be the space of bounded continuous functions.

Since the choice of $T>0$ in Theorem \ref{Th:1} is arbitrary, the solution exists for all $t \geq 0$, thus $y(t)$ also exists for all $t \geq 0$. Replacing the initial interval $[-h,0]$ with $[-h+s,s]$ for all $s \geq  0$, we can guarantee the existence and uniqueness of the solutions for $t \geq s \geq 0$ with the initial $\mathcal{F}_s$-measurable functions $\f(\theta, \cdot), \f(0, \cdot)$, which satisfy the conditions of Theorem \ref{Th:1} on $[s-h,s]$. This solution will be denoted with $u(t,s,\f)$. Similarly,
\[
u_t(s,\f) = u(t + \theta, s, \f), \theta \in [-h,0]
\]
is a shift of the solution $u(t,\f)$, such that $u_s(s,\f) = u(s+\theta, s, \f) = \f(\theta)$ and for $\theta = 0$, $\f(0,\cdot) \in B_0^{\rho}$.

Following \cite{Car}, define the family of shift operators
\begin{equation}\label{family}
U_{s}^t \f:= u(t+\theta,s,\f) = u_t(s,\f).
\end{equation}

Denote $\F_{s}^{t}(d W)$ to be the minimal $\sigma$-algebra  containing $W(\tau) - W(s), \tau \in [s,t]$. Note that $u_t(s,\f)$ is independent of the $\sigma$-algebra $G^t$, which is the minimal sigma-algebra containing $W(\tau) - W(t)$ for $\tau \geq t$.

For any nonrandom $\f \in B^{\rho}$ with $s \geq 0$ and $t \geq s$, $U_{s}^{t} \f:= u_t(s,\f)$ is an $\mathcal{F}_{s}^{t} (d W)$ measurable random function taking values in $B_1^{\rho}$, with $u(t,s,\f) \in B_0^{\rho}$ for $\theta = 0$. Defining
$y(t,s,\f) = (u(s,t,\f), u_t(s,\f))$, we have that $y$ maps $B^{\rho}$ into itself.

The next proposition follows from Theorem \ref{Th:1}:
\begin{proposition}\label{Prop:2}
The family of the operators (\ref{family}) satisfies
\begin{equation}\label{dyn_sys}
    U_{\tau}^{t} U_s^{\tau} \f = U_s^t
\end{equation}
for all $t \geq \tau \geq s \geq 0$ and $\f \in B^{\rho}$.
\end{proposition}
Let $D$ be a $\sigma$-algebra of Borel subset of $B^{\rho}$. Then $y(t,s,\f)$ naturally denotes the following probability measure $\mu_t$ on $D$:
\begin{equation}\label{measure}
    \mu_t(A)  = P\{y(t,s,\f) \in A\} = P\{U_s^t \f \in A\} = P(s,\f,t,A)
\end{equation}
The measure $\mu$ is the transition function corresponding to the random process $y(t,s,\f)$. In a similar way as in the finite dimensional case \cite{Car}, one can show that this function satisfies the properties of the transition probability. This way we have

\begin{theorem}(Markov property)
Under the assumptions of Theorem \ref{Th:1}, the process $y(t,s,\f) \in B^{\rho}$ is the Markov process on $B^{\rho}$ with the transition function $P(s,\f,t,A)$ given by (\ref{measure}).
\end{theorem}

\begin{proposition}\label{Prop:3}
For any $t \geq s \geq 0$ we have
\[
P(s,\f, t, A) = P(0, \f, t-s, A)
\]
\end{proposition}
\begin{proof}
Let $\tilde{u}(t) = u(s+t,s,\f)$. Then $\tilde{u}(0) = \f(0,\cdot)$ and $\tilde{u}_0 = u(s+\theta, s, \f) = \f(\theta, \cdot)$. On the other hand,
\[
\tilde{u}(t) = u(s+t,s,\f) = S(t) \f(0,\cdot) +
\int_{s}^{s+t} S(s+t - \tau) f(u(\tau)) d \tau +
\]
\[
+\int_{s}^{s+t} S(s+t - \tau) \sigma(u(\tau)) d W(s+\tau) = S(t) \f(0,\cdot) +
\]
\[
+ \int_{0}^{t} S(t - \tau) f(u(\tau+s)) d \tau + \int_{0}^{t} S(t - \tau) \sigma (u(\tau+s)) d \tilde{W}(\tau)
\]
where
\[
\tilde{W}(\tau) := W(s+\tau) - W(s)
\]
is once again a $Q$-Wiener process. This way $\tilde{u}$ solves
\begin{equation}\label{eqn0}
\tilde{u}(t) = S(t) \f(0,\cdot) +
\int_{0}^{t} S(t - \tau) f(\tilde{u}(\tau)) d \tau + \int_{0}^{t} S(t - \tau) \sigma (\tilde{u}(\tau)) d \tilde{W}(\tau)
\end{equation}
The same equation is satisfied with $u(t,0,\f)$  such that $u(0,0,\f) = \f(0,\cdot)$ and $u_0 = \f(\theta,\cdot)$. The only difference is that $u(t,0,\f)$ solves (\ref{eqn0}) with a different Wiener process $W$. However, since the distribution of $W$ is the same as $\tilde{W}$, the distribution of $u(s+t,s,\f)$ is the same as the distribution of $u(t,0,\f)$, and hence independent of $s$. Thus the distribution of $u_t(s,\f) = u(t+ \theta, s,\f) = u(t-s + \theta +s, s \f)$ coincides with the distribution of $u(t-s+\theta, 0, \f) = u_{t-s}(0,\f)$. Hence
\[
P(s,\f,t,A) = P\{u_t(s,\f) \in A\} = P\{u(t+\theta,s,\f) \in A\} =
\]
\[
=P\{u(t-s+\theta,0,\f)\in A\} = P\{u_{t-s}(0,\f) \in A\}
\]
which yields the desired result.
\end{proof}
For $g \in B_b(B^{\rho})$, for $\f \in B^{\rho}$ and $t \geq s \geq 0$ define
\[
P_{s,t}(\f) := \E g(y(t,s,\f)).
\]
From proposition \ref{Prop:3} we have  $P_{0,t-s}(\f)$ and denote $P_t \f = P_{0,t}(\f)$.
From Theorem \ref{Th:2} and Proposition \ref{Prop:1} we have
\begin{proposition}\label{Prop:4}
Under the assumptions of Theorem \ref{Th:1} the transition semigroup $P_t, t \geq 0$ is stochastic continuous an satisfies the Feller property
\[
P_t: C_b(B^\rho) \to C_b(B^\rho) \text{ and } \lim_{t \to 0} P_t \f(\theta) = \f(\theta).
\]
\end{proposition}
Define $\bar{\rho}(x) = (1+|x|^{\bar{r}})^{-1}$. The main result of the paper is the following theorem:
\begin{theorem}\label{Th:4IM}
Let the assumptions of Theorem \ref{Th:1} hold. Assume the equation (\ref{MildSOl}) has a solution in $B^{\bar{\rho}}$ which is bounded in probability for $t \geq 0$ with $r>d+\bar{r}$. Then there exists an invariant measure $\mu$ on $B^{\rho}$, i.e.
\[
\int_{B^{\rho}} P_t \f(x) d \mu(x) = \int_{B^{\rho}} \f(x) d\mu, \text{ for any } t \geq 0 \text{ and } \f \in C_b(B^{\rho}).
\]
\end{theorem}

The key condition in Theorem \ref{Th:4IM} is the existence of a globally bounded solution. The next theorem provides the sufficient conditions for the existence of such solution in terms of the coefficients, in the case when $A$ is in the divergence form.

\begin{theorem}\label{Th:5}
Assume
\begin{itemize}
\item $D = \R^d, d \geq 3$;
\item the conditions of Theorem \ref{Th:1} hold;
\item for some $\sigma_0 > 0$, we have $|\sigma(u)| \leq \sigma_0$, $\forall u \in B_1^{\rho}$;
\item there exists $\Psi \in L^1(\R^d) \cap L^{\infty}(\R^d)$ such that $|f(u(x))| \leq \Psi(x)$ for all $u \in B_1^{\rho}, x \in \R^d$;
\item $u(t,\cdot) = \f(t,\cdot), t \in [-h,0]$, $u(0,\cdot) = \f(0,x)$ satisfy
\[
\E \int_{\R^d} |\f(0,x)|^2 dx < \infty \text{ and } \E \int_{\R^d} \int_{-h}^0|\f(\theta,x)|^2 dx d\theta < \infty.
\]
\end{itemize}
Then
\[
\sup_{t \geq 0} \E \|y(t)\|^2_{B^{\rho}} < \infty,
\]
which is a sufficient condition  for the boundedness in probability.
\end{theorem}

Finally, in the case when the domain $D$ is bounded, we may establish the uniqueness of the stationary solution as well as its stability. In this section, the weight $\rho \equiv 1$, thus
\[
B_0:= L^2(D)
\]
\[
B_1:=L^2(-h,0, B_0)
\]
and
\[
B:=B_0 \times B_1.
\]
The semigroup (\ref{def_semi}) now satisfies the exponential estimate
\[
\|S(t) u_0\|_{B_0}^2 \leq e^{-2 \lambda t}\|u_0\|_{B_0}^2,
\]
where $\lambda_1>0$ is the principle eigenvalue of $-A$. 
In a standard way, we may extend the Q-Weiner process $W(t)$ to $t \in \R$ as 
\[
W(t) = \begin{cases}
W(t), t \geq 0;\\
V(-t), t \leq 0.
\end{cases}
\]
Here $V$ is another Q-Weiner process, independent of $W$.  
\begin{definition}
A $B_0$-valued process $u(t)$ is a mild solution of (\ref{MainB}) for $t \in \R$ if 
\begin{enumerate}
\item for all $t \in \R$, $u_t$ is $\mathcal{F}_t$ measurable;
\item for all $t \in \R$ 
\[
\E \|u(t)\|^2_{B_0} < \infty;
\]
\item for all $-\infty < t_0 < t < \infty$ with probability 1 we have 
\[
u(t) = S(t-t_0) u(t_0) +\int_{t_0}^t S(t-s) f(u_s) ds + \int_{t_0}^t S(t-s) \sigma(u_s) d W(s)
\]
\end{enumerate}
\end{definition}
\begin{theorem}\label{Th:6}
Assume the Lipschitz constant $L$ is sufficiently small (see (\ref{smallness}) for the exact condition), then the equation (\ref{MainB}) has a unique solution $u^*(t,x)$, defined for $t \in \R$, and 
\[
\sup_{t \in \R} \E\|u^*(t)\|^2_B< \infty. 
\]
Furthermore, this solution is exponentially attractive, that is exist $K, \gamma >0$ such that for all $t_0 \in \R$ and $t> t_0 + h$, and for any other solution $\eta(t)$ with $\eta(t_0) \in B_0$ and $\eta_{t_0} \in B_1$ we have
\[
\E \|u(\cdot,t) - \eta(\cdot,t)\|_B^2 \leq K e^{-\gamma(t-t_0)} \E \|u(\cdot,t_0) - \eta(\cdot,t_0)\|_B^2.
\]
\end{theorem}
\section{Proofs of the Main Results}
\subsection{Proof of Theorem \ref{Th:1}} Let
$B_{p,T}$, $p \geq 2$ be the space of $\mathcal{F}_t$-measurable for $t \in [0,T]$ processes, equipped with the norm $\|\Phi\|_{B_{p,T}} := \E \int_{-h}^{T} \|\Phi(t,\cdot)\|_{B_0^{\rho}}^p dt$, $\Phi: B_0^{\rho} \to B_0^{\rho}$. Define
\begin{equation}\label{Psi}
    \Psi \Phi(t,\cdot):= S(t) \Phi(0,\cdot) + \int_0^t S(t-s)f(\Phi(s+\theta, \cdot)) ds + \int_0^t S(t-s) \sigma(\Phi(s+\theta, \cdot)) dW(s)
\end{equation}
for $t \in [0,T]$, and
\[
\Psi \Phi(t,\cdot) = \f(t,\cdot), t \in [-h,0]; \text{ with } \Psi \Phi(0,\cdot) = \f(0,\cdot).
\]
This way
\[
\Psi \Phi(t,\cdot)\|^p_{B_{p,T}} \leq \E \int_{-h}^0\|\f(t,\cdot)\|^p_{B_0^{\rho}} dt + 3^{p-1} \E \int_0^T \|S(t) \f(0,\cdot)\|_{B_0^{\rho}}^p dt +
\]
\[
+3^{p-1} \E \int_0^T \left\|\int_0^t S(t-s)f(\Phi(s+\theta, \cdot)) ds \right\|^p_{B_0^{\rho}} dt + 3^{p-1} \E \int_0^T \left\|\int_0^t S(t-s)\sigma(\Phi(s+\theta, \cdot)) d W(s) \right\|^p_{B_0^{\rho}} dt \leq
\]
\[
\leq C_1(T) + 3^{p-1} (I_1+ I_2 + I_3).
\]
It follows from (\ref{ast}) that
\[
I_1 \leq C_{\rho}^p(T)  \int_0^T \E \|\f(0,\cdot)\|^p_{B_0^{\rho}} dt < \infty;
\]
Next, using the conditions [i] and [ii] for $f$, we have
\[
I_2 \leq C_{\rho}^p(T) \int_0^T T^{p-1} \left(\E  \int_0^t \|f(\Phi_s,\cdot)\|_{B_0^{\rho}}^p ds \right) dt \leq C_2 \int_0^T dt \int_0^t\left(1+\E \|\Phi_s\|^p_{B_1^{\rho}}\right) ds \leq
\]
\begin{equation}\label{I2}
\leq C_3+C_2 \int_{0}^{T} dt \left(\int_0^t \E \left(\int_{-h}^0 \|\Phi(s+ \theta, \cdot) \|^2_{B_0^{\rho}}\|\right)^{\frac{p}{2}}\right) \leq C_3 + C_4 \E \int_{-h}^T \|\Phi(t,\cdot) \|_{B_0^{\rho}}^p dt < \infty.
\end{equation}
In order to estimate $I_3$, we use Lemma 7.2 \cite{DapZab92} and (\ref{Noise_est}). We have
\[
I_3 \leq C(P) \int_0^{T} \E \left( \int_0^t\|S(t-s) \sigma(\Phi_s(\cdot))\|^2_{\mathcal{L}_2} ds \right)^{\frac{p}{2}}dt
\]
\begin{equation}\label{I3}
\leq C(p) a^p \sup_n \|e_n\|_{\infty}^p  \int_0^T \E \left( \int_0^t\|S(t-s) \sigma(\Phi_s(\cdot))\|^2_{B_0^\rho} ds \right)^{\frac{p}{2}} dt  \leq C_4 + C_5 \int_0^T \int_0^t \E \|\Phi_s(\cdot) \|_{B_1^{\rho}}^p ds < \infty
\end{equation}
the same way as in (\ref{I2}). Combining these estimates, we have $\Psi: B_{p,T} \to B_{p,T}$.
We next show that $\Psi$ is contractive.
For any $\Phi, \tilde{\Phi} \in B_{p,t}$ we
\[
\|\Psi \Phi(s,\cdot) - \Psi \tilde{\Phi}(s,\cdot)\|_{B_{p,T}}^p \leq 2^{p-1} \int_0^t \E  \left\|\int_0^s S(s-\tau)(f(\Phi_\tau(\cdot) - f(\tilde{\Phi}_\tau(\cdot)) d\tau \right\|_{B_0^{\rho}}^p ds +
\]

\begin{equation}\label{Stisk}
+2^{p-1} \int_0^t \E  \left\|\int_0^s S(s-\tau)(\sigma(\Phi_\tau(\cdot)) - \sigma(\tilde{\Phi}_\tau(\cdot)) d\tau \right\|_{B_0^{\rho}}^p ds := 2^{p-1} (I_4 + I_5)
\end{equation}
\[
I_4 \leq C_\rho^p(T) L^p \int_0^t \E \left(\int_0^s\|\Phi_\tau(\cdot) - \tilde{\Phi}_\tau( \cdot)\|_{B_1^\rho} d \tau \right)^p ds \leq
\]
\begin{equation}\label{I4}
\leq C(\rho, T,\rho) \int_0^t \int_0^s \E \left( \int_{-h}^0 \|\Phi(\tau+ \theta, \cdot) - \tilde{\Phi}(\tau+ \theta,\cdot)\|^2_{B_0^{\rho}} d\theta  \right)^{p/2}  d\tau ds \leq C_5(\rho, T,p) t^2 \|\Phi-\tilde{\Phi}\|^p_{B_{p,t}}
\end{equation}
Now following the estimate (\ref{I3}), we have
\[
I_5  \leq C(P) a^p \sup_n\|e_n\|_{\infty}^p \int_0^t \E  \left(\int_0^s \| S(s-\tau) [\sigma(\Phi_\tau) - \sigma(\tilde{\Phi}_\tau)]\|^2_{B_0^{\rho}}\right)^{\frac{p}{2}} dt \leq
\]

\begin{equation}\label{I5}
\leq C_6  \int_0^t \left( \int_0^s \E \left(\int_{-h}^0 \|\Phi(\tau + \theta, \cdot) - \tilde{\Phi}(\tau + \theta, \cdot) \|_{B_0^{\rho}}^2 d\theta\right)^{p/2} d \tau \right) ds \leq C_6(\rho,T,p,h) t^2 \|\Phi-\tilde{\Phi}\|^p_{B_{p,t}}.
\end{equation}

Consequently, for $\tilde{t}$ small enough, (\ref{I4}) and (\ref{I5}) imply that the map $\Psi$ has a unique fixed point in $B_{p,\tilde{t}}$, which is the solution of \eqref{MildSOl}. If we consider the problem on $[0,\tilde{t}], [\tilde{t}, 2\tilde{t}],...$ with $C_6 \tilde{t}^2<1$. Since the solution is continuous with probability 1 in $B_0^{\rho}$ norm, we get the existence and uniqueness of the solution on $[0,T]$.

It remains to prove the estimate (\ref{cont dep}). It follows from (\ref{MildSOl}) that for any $t \in [-h,T]$ we have
\[
\E \|u(t,\cdot)\|_{B_0^{\rho}}^p \leq 3^{p-1} \E \|S(t) \f(0,\cdot)\|^p_{B_0^\rho} +
\]
\[
+ 3^{p-1} \E \left(\int_0^t \|S(t-s) f(u(s))\|_{B_0^{\rho}} ds \right)^p
+ 3^{p-1} \E \left\|\int_0^t S(t-s) \sigma(u(s)) d W(s) \right\|_{B_0^{\rho}}^p \leq
\]
\[
\leq 3^{p-1} C_{\rho}(T) \E \|\f(0,\cdot)\|_{B_0^\rho}^p + 3^{p-1} C_7  \int_0^t(1+ \E \|u(s)\|^p_{B_1^{\rho}}) ds + 3^{p-1} C_8 \E \left(\int_0^t \|S(t-s) \sigma(u(s))\|^2_{\mathcal{L}^2} ds\right)^{\frac{p}{2}} 
\]
\begin{equation}\label{3.9}
\leq C_9 \left(\E \|\f(0,\cdot)\|^p_{B_0^{\rho}} +  \int_0^t(1+ \E \|u(s)\|_{B_1^\rho}^p) ds\right).
\end{equation}
Consider two separate cases $t \in [0,h]$ and $t \in [h,T]$:\\
If $t \in [0,h]$,  then
\[
\E \|u(t)\|^p_{B_1^\rho} = \E \left(\int_{-h}^0 \|u(t+ \theta, \cdot)\|^2_{B_0^{\rho}} d \theta \right)^{\frac{p}{2}} \leq 2^{\frac{p}{2} -1} \left(\E \left(\int_{-h}^{-t} \|u(s, \cdot)\|_{B_0^{\rho}}^2 ds  \right)^{\frac{p}{2}} + \E \left(\int_0^t \|u(s,\cdot)\|_{B_0^\rho}^2 ds\right)^\frac{p}{2}\right)
\]
\begin{equation}\label{3.10}
\leq 2^{\frac{p}{2}-1}\left(\E  \|\f(t,\cdot)\|^p_{B_1^\rho} + h^{\frac{p-2}{p}} \int_0^h \E \|u(s,\cdot)\|_{B_0^\rho}^p ds\right)
\leq
2^{\frac{p}{2}-1}\left(\E  \|\f(t,\cdot)\|^p_{B_1^\rho} + C_{10}\sup_{s \in [0,t]} \E \|u(s,\cdot)\|_{B_0^\rho}^p\right).
\end{equation}
If $t \in [h, T]$, we have
\begin{equation}\label{3.11}
    \E  \|u(t)\|_{B_1^{\rho}}^p = \E \left(\int_{-h}^0 \|u(t+ \theta, \cdot) \|^2_{B_0^\rho} d \theta\right)^\frac{p}{2} \leq C_{11}(T) \sup_{s \in [0,t]} \E \|u(s)\|_{B_0^{\rho}}^p.
\end{equation}
From (\ref{3.9}-\ref{3.11}) we have
\[
\sup_{s \in [0,t]}\E \|u(s,\cdot)\|_{B_0^{\rho}}^p \leq C_{12}(T) \left(\E \|\f(0,\cdot)\|^p_{B_0^{\rho}} + \E \|\f(t,\cdot)\|^p_{B_1^{\rho}} + \int_0^t \sup_{\tau \in [0,s]} \E \|u(\tau, \cdot) \|^p_{B_0^{\rho}} ds\right).
\]
Estimating the last term separately, we have
\begin{equation*}\label{3.12}
\sup_{s\in [0,t]} \E \|u(s,\cdot)\|_{B_0^{\rho}}^p \leq C_{13}(T)[1+ \E\|\f(0,\cdot)\|_{B_0^{\rho}}^p +  \E \|\f(t, \cdot)\|_{B_1^\rho}^p].
\end{equation*}
Combining the estimates above, we get
\[
\E \|u(t)\|_{B_1^\rho}^p \leq C_{14}(T)(1+\E \|y(0)\|^p_{B^\rho},
\]
which completes the proof.
\subsection{Proof of Theorem \ref{Th:2}} By definition of $y_1$ and $y_2$,
\begin{equation}\label{3.13}
    \sup_{t \in [0,T]} \E\|y(t) - y_1(t)\|_{B^\rho}^2 \leq \sup_{t\in [0,T]} \E \|u(t,\phi) - u(t,\phi_1)\|_{B_0^{\rho}}^2 + \sup_{t \in [0,T]} \E \|u_t(\phi) - u_t(\phi_1)\|^2_{B_1^\rho}.
\end{equation}
The first term in (\ref{3.13}) can be estimated as follows
\begin{equation}\label{3.14}
    \sup_{t \in [0,T]} \E \|u(t,\phi) - u(t,\phi_1)\|_{B_0^\rho}^2 \leq C_{15} \sup_{t \in [0,T]} \E \|\phi(t) - \phi_1(t)\|^2_{B_0^\rho}.
\end{equation}
As for the second term in (\ref{3.13}), once again we consider separately the cases $t \in [0,h]$ and $t \in [h,T]$, and taking into account the estimate (\ref{3.14}), we get
\[
\sup_{t \in [0,T]} \E \int_{-h}^{0}\|u(t+\theta,\phi) - u(t+\theta, \phi_1)\|^2_{B_0^\rho} d \theta \leq C_{16} \sup_{t \in [0,T]} \E \|\phi(t) - \phi_1(t)\|_{B^\rho}^2,
\]
which completes the proof.

\subsection{Proof of Theorem \ref{Th:4IM}} We start with the following auxiliary lemmas.
\begin{lemma}\label{L3}
For any fixed $T_0 > 2h$, the operator
\[
A \f_0:= S(T_0+ \theta): B_0^{\bar{\rho}} \to B_1^\rho
\]
is a Hilbert-Schmidt operator.
\end{lemma}
\begin{proof}
By \cite{OvsPel}, there exists an orthonormal basis $\{h_n, n \geq 1\}$ in $B_0^\rho$ such that  $\sup_n \|h_n\|_{L^\infty(D)} < \infty$. It is straightforward to verify that if $\{e_n,n \geq 1\}$ is an orthonormal basis in $H = L^2(D)$, then $\{\frac{e_n}{\bar{\rho}^{1/2}},  n \geq 1\}$ is an orthonormal basis in $B_0^{\bar{\rho}}$. Therefore

\[
\|A\|^2_{\mathcal{L}^2} = \sum_{i=1}^\infty \left\|A  \frac{e_i}{\sqrt{\bar{\rho}}} \right\|^2_{B_1^\rho} = \sum_{i=1}^\infty \left\|S(T_0 + \theta) \frac{e_i}{\sqrt{\bar{\rho}}}\right\|^2_{B_1^\rho} = \sum_{i=1}^\infty \int_{-h}^0 d\theta \int_D \left|S(T_0 + \theta) \frac{e_i}{\sqrt{\bar{\rho}}}\right|^2 \rho(x) dx  =
\]

\[
= \sum_{i=1}^\infty \int_{-h}^0 d\theta \int_D \left|\int_D G(T_0+ \theta,x,y) \frac{e_i(y)}{\sqrt{\bar{\rho}(y)}} dy\right|^2 \rho(x) dx = \int_{-h}^0 d\theta \int_D \rho(x) \int_{D}\frac{G^2(T_0+\theta,x,y)}{\bar{\rho}(y)} dy dx \leq
\]

\[
\leq \int_{-h}^0 d \theta \int_D \frac{\rho(x)}{\bar{\rho}(y)} \int_{D} C_1(T_0) (T_0+ \theta)^{-d} \exp\{-2C_2(T_0) \frac{|x-y|^2}{T_0 + \theta}\} dy dx \leq
\]

\[
\leq C_{17} \int_{-h}^0 \frac{d \theta}{(T_0 + \theta)^{d/2}} \int_{\R^d}  \int_{\R^d} \frac{1}{(T_0+ \theta)^{d/2}}\exp\{-2C_2(T_0) \frac{|x-y|^2}{T_0 + \theta}\} \frac{\rho(x)}{\bar{\rho}(y)} dx dy.
\]
But

\[
\int_{\R^d} \frac{1}{(T_0+ \theta)^{d/2}}\exp\{-2C_2(T_0) \frac{|x-y|^2}{T_0 + \theta}\} \frac{\rho(x)}{\rho(y)} dx \rho(y) \leq
\]
\[
\leq C(r)\int_{\R^d} \frac{1}{(T_0+ \theta)^{d/2}}\exp\{-2C_2(T_0) \frac{|x-y|^2}{T_0 + \theta}\} (1+|x-y|^r) dx \rho(y) \leq C_{18}(T,r) \rho(y).
\]
Thus
\[
\|A\|^2_{\mathcal{L}^2} \leq C_{19}(T,r) \int_{-h}^0 \frac{d \theta}{(T_0 + \theta)^{d/2}} \int_{\R^d} \frac{1+|y|^{\bar{r}}}{1+|y|^r} dy < \infty,
\]
which completes the proof of the Lemma.
\end{proof}
\begin{corollary}\label{cor1}
Following the lines of the proof of Lemma \ref{L3} we can show that $S(t)$ is a compact operator from $B_0^{\bar{\rho}}$ to $B_0^\rho$ for $t>0$.
\end{corollary}
We now return to the proof of Theorem \ref{Th:4IM}. We will follow the approach of the Theorem 11.29 from \cite{DapZab92}. We have
\begin{equation}\label{3.16}
    u(T_0) = S(T_0) \f(0, \cdot) + \int_0^{T_0} S(T_0 - s) f(u_s) ds + \int_{0}^{T_0} S(T_0 - s) \sigma(u_s) dW(s)
\end{equation}
\begin{equation}\label{3.17}
    u_{T_0} = u(T_0 + \theta) = S(T_0+ \theta) \f(0, \cdot) + \int_0^{T_0+\theta} S(T_0 + \theta - s) f(u_s) ds + \int_{0}^{T_0+ \theta} S(T_0 + \theta - s) \sigma(u_s) dW(s)
\end{equation}
The arguments in Theorem 11.29 from \cite{DapZab92} can be applied to (\ref{3.16}) directly.
\begin{lemma}\label{L4}
For $p>2$ and $\alpha \in \left(\frac{1}{p}, \frac{1}{2}\right)$, the operator
\[
(G_{\alpha} \f)(\theta) = \int_0^{T_0+ \theta} (T_0 + \theta -s)^{\alpha -1} S(T_0 + \theta -s) \f(s) ds
\]
is a compact operator from $L^p(0,T_0, B_0^{\bar{\rho}})$ into $C(-h, 0, B_0^\rho)$.
\end{lemma}
\begin{remark}
Compactness in $C(-h,0, B_0^\rho)$ implies compactness in $B_1^\rho$.
\end{remark}
\begin{proof}
We will make use of the infinite dimensional version of Arzela-Ascoli Theorem. To this, we need to show
\begin{itemize}
    \item[{[i]}] For any fixed $\theta \in [-h,0]$ the set
    \[
    \{G_{\alpha}(\f)(\theta), \|\f\|_{L^p} \leq 1\}
    \]
     is compact in $B_0^\rho$;
    \item[{[ii]}] for any $\varepsilon > 0$ there exists $\delta>0$ such that if $\|\f\|_{L^p} \leq 1$ such that if $\|\f\|_{L^p} \leq 1$ and $\forall \theta_1, \theta_2$ with $|\theta_1 - \theta_2| \leq \delta$ we have
\[
    \|G_{\alpha}(\f)(\theta_1) - G_{\alpha}(\f)(\theta_2)\|_{B_0^\rho} < \ve.
    \]
\end{itemize}
To check [i], for fixed $\theta \in [-h,0]$ and $\ve>0$ introduce
    \[
    G_{\alpha}^\ve \f := \int_0^{T_0 + \theta -\ve} (T_0 + \theta -s)^{\alpha -1} S(T_0 + \theta -s) \f(s) ds = S(\ve) \int_0^{T_0 + \theta -\ve} (T_0 + \theta -s)^{\alpha -1} S(T_0 + \theta -s - \ve) \f(s) ds
    \]
    Clearly
    $\int_0^{T_0 + \theta -\ve} (T_0 + \theta -s)^{\alpha -1} S(T_0 + \theta -s - \ve) \f(s) ds$ is in $B_0^{\bar{\rho}}$. Using Corollary \ref{cor1}, $S(\ve)$ is a compact operator from $B_0^{\bar{\rho}}$ to $B_0^\rho$. Then, following \cite{DapZab92}, p.227, $G_\alpha^\ve$ converges to $G_\alpha$ strongly as $\ve \to 0$, hence $G_\alpha$ is compact hence [i] follows.\\
    To prove [ii], fix $\theta$ and $r$ such that $-h \leq \theta \leq \theta + r \leq 0$, and $\|\f\|_{L^p}\leq 1$. Then
    \[
    \|(G_\alpha \f)(\theta + r) - (G_\alpha \f)(\theta)\|_{B_0^\rho}
    \]
    \[
     =\left\|\int_0^{T_0+ \theta + r} (T_0 + \theta + r-s)^{\alpha -1} S(T_0 + \theta + r -s ) \f(s) ds - \int_0^{T_0+ \theta} (T_0 + \theta -s)^{\alpha -1} S(T_0 + \theta -s ) \f(s) ds\right\|_{B_0^\rho}
    \]
    \[
    \leq \int_0^{T_0+ \theta} \left\|(T_0 + \theta + r-s)^{(\alpha -1)} S(T_0 + \theta + r -s ) - (T_0 + \theta -s)^{(\alpha -1)} S(T_0 + \theta -s)\right\| \|\f(s)\| ds +
    \]
    \[
    + \int_{T_0 + \theta}^{T_0 + \theta + r} \left\|(T_0 + \theta + r-s)^{(\alpha -1)} S(T_0 + \theta + r -s ) \f(s)\right\| ds \leq
    \]
    \[
    \leq \left(\int_0^{T_0} \|(r+s)^{\alpha-1} S(s+r) - s^{\alpha -1} S(s) \|^q ds\right)^{\frac{1}{q}} \|\f\|_{L^p} + C_{20} \left(\int_0^{T_0}s^{(\alpha -1)q} ds\right)^{\frac{1}{q}} \|\f\|_{L^p}:=J_1 + J_2.
    \]
    Direct calculation yields
    \[
    J_2 = C_{20} \frac{r^{\alpha - \frac{1}{p}}}{((\alpha-1)q+1)^{\frac{1}{q}}} \|\f\|_{L^p} \to 0 \text{ as } r\to 0.
    \]
    We now proceed with estimating $J_1$. Since $S(t)$ is compact, then $S(t)$ is strongly continuous for $t>0$ (see \cite{Paz}, Theorem 3.27), hence $\|S(s+r)-S(s)\| \to 0$ as $r \to 0$, for any $s>0$. Furthermore, the integrand in $J_1$ is bounded by $2C_{20} s^{(\alpha-1)q}$. Hence, by Dominated Convergence Theorem, $J_1 \to 0$ as $r \to 0$, which concludes the proof of the Lemma.

\end{proof}
We may now complete the proof of Theorem \ref{Th:4IM}. For any $r>0$ introduce
\[
K(r):=\{(\mu, \nu),  \mu \in B_0^\rho, \nu \in B_1^\rho\}
\]
such that
\[
\mu := S(T_0) v  + (G_1 \f)(0) + (G_{\alpha} h)(0);
\]
\[
\nu := S(T_0 + \theta) v + (G_1 \f)(0) + (G_{\alpha} h)(0)
\]
with $\|v\|_{B_0^\rho} \leq r$, $\|\f\|_{L^p(0,T_0,B_0^{\bar{\rho}})} \leq r$ and $\|h\|_{L^p(0,T_0, B_0^{\bar{\rho}})} \leq r$.
It follows from Lemma \ref{L3}, Corollary \ref{cor1} and Lemma \ref{L4} that $K(r)$ is compact in $B^\rho$.
\begin{lemma}\label{L5}
Under the conditions of Theorem \ref{Th:1}, there is $C>0$ such that for arbitrary $r>0$ and $y = (x,z) \in B^{\bar{\rho}}$ such that
$\|y\|_{B^{\bar{\rho}}} \leq r$ we have
\begin{equation}\label{3.19}
    P\{(u(T_0,x,z), u_{T_0}(x,z)) \in K(r)\} \geq 1 - c r^{-p}(1+ \|y\|^p_{B^{\bar{\rho}}})
\end{equation}
\end{lemma}
\begin{proof}
From factorization formula \cite{DapZab96}, Th. 5.2.5, we have
\begin{equation}\label{3.20}
u(T_0,y) = S(T_0)x + (G_1 f(u_s))(0) +\frac{\sin(\alpha \pi)}{\pi}(G_\alpha Y(s))(0),
\end{equation}
\begin{equation}\label{3.21}
u_{T_0}(y) = S(T_0+\theta)x + (G_1 f(u_s))(\theta) + \frac{\sin(\alpha \pi)}{\pi}(G_\alpha Y(s))(\theta)
\end{equation}
and
\begin{equation*}
    Y(s) = \int_0^s (s-\tau)^{-\alpha} S(s-\tau) \sigma(u_\tau) d W(\tau)
\end{equation*}
Using Lemma 7.2 \cite{DapZab92}, we obtain
\[
\E \int_0^{T_0} \|Y(s)\|^p_{B_0^{\bar{\rho}}} ds = \E \int_0^{T_0} \left\|\int_0^s (s-\tau)^{-\alpha} S(s-\tau) \sigma(u_\tau) d W(\tau)\right\|_{B_0^{\bar{\rho}}}^p ds \leq
\]
\[
\leq C_{p,T_0} \E \int_0^{T_0}\left(\int_0^s(s-\tau)^{-2\alpha}\left\|S(s-\tau) \sigma(u_\tau) \circ Q^{1/2}\right\|_{\mathcal{L}_2(H, B_0^{\bar{\rho}}}\right)^{\frac{p}{2}} ds\leq
\]
\begin{equation}\label{3.22}
\leq C_{21} \E \int_0^{T_0}\left(\int_0^s(s-\tau)^{-2\alpha}\|\sigma(u_\tau)\|^2_{B_0^{\bar{\rho}}} d\tau\right)^{\frac{p}{2}}.
\end{equation}
Using Hausdorff-Young's inequality and \eqref{MildSOl}, we have
\begin{multline}\label{3.23}
\E \int_0^{T_0} \|Y(s)\|^p_{B_0^{\bar{\rho}}} \leq C_{21} \left(\int_0^{T_0} t^{-2 \alpha} dt\right)^{\frac{p}{2}} \int_{0}^{T_0} \E \|\sigma(u_t)\|^p_{B_0^{\bar{\rho}}} dt \leq C_{22} \int_0^{T_0}(1+ \E\|u_t\|^p_{B_1^{\bar{\rho}}}) dt \\
\leq C_{23}(1+ \|y\|^p_{B^{\bar{\rho}}}).
\end{multline}
In a similar way,
\begin{equation}\label{3.24}
    \E  \int_0^{T_0} \|f(u_s)\|^p_{B_0^{\bar{\rho}}} ds \leq C_{23}(1+ \|y\|^p_{B^{\bar{\rho}}}).
\end{equation}
Hence, if $\|y\|_{B^{\bar{\rho}}} \leq r$, $\|f(u_s)\|_{L^p(0,T_0,B_0^{\bar{\rho}})} \leq r$ and 
$$\|\sigma(u_s)\|_{L^p(0,T_0,B_0^{\bar{\rho}})} \leq \frac{\pi r}{\sin(\alpha \pi)}$$
then from the definition of $K(r)$ we have
\[
(u(T_0,y), u_{T_0}(y)) \in K(r).
\]
Assume $\|y\|_{B^{\bar{\rho}}} \leq r$. Then
\[
P\{(u(T_0,y), u_{T_0}(y)) \notin K(r)\} \leq P\{\|f(u_s)\|_{L^p(0,T_0; B_0^{\bar{\rho}})} > r\} + P\{ \|y(s)\|_{L^p(0,T_0,B_0^{\bar{\rho}})}\} \leq 2 r^p C_{23}(1+\|y\|^p_{B^{\bar{\rho}}})
\]
where we used (\ref{3.23}) and (\ref{3.24}). Hence the proof of Lemma \ref{L5} follows.
\end{proof}
The rest of the proof of Theorem \ref{Th:4IM} follows the lines of the proof of Theorem 11.29, \cite{DapZab92}.

\subsection{Proof of Theorem \ref{Th:5}} The proof of this theorem has a lot in common with Theorem 1, \cite{MisStaYip}. Let us point out the differences due to the presence of the delay. We have
\begin{equation}\label{3.26}
\E  \|y(t)\|^2_{B^\rho} = \E \int_{D} |u(t,x)|^2 \rho(x) dx + \E \int_{-h}^0 d \theta \int_D |u(t+\theta,x)|^2 \rho(x) dx.
\end{equation}
By definition of a mild solution \eqref{MildSOl}, we have
\[
\|u(t,x)\|^2_{B_0^\rho} \leq 3(I_1(t) + I_2(t) + I_3(t))
\]
where
\[
I_1(t) = \int_{\R^d} \left(\int_{\R^d} G(t,x,y) \f(0,y) dy \right)^2 \rho(x) dx,
\]
\[
I_2(t) = \int_{\R^d} \left(\int_0^t  \int_{\R^d} G(t-s,x,y) f(u_s(y)) dy ds\right)^2 \rho(x) dx,
\]
\[
I_3(t) = \int_{\R^d} \left(\int_{\R^d} G(t-s,x,y) \sigma(u_s(y)) d W(s) dy\right)^2 \rho(x) dx.
\]
It follows from (\ref{GFest}) that for all $t \geq 0$
\[
\E I_1 \leq \int_{\R^d} \left(\int_{\R^d} G(t,x,y) dy \int_{\R^d} G(t,x,y) dy \f^2(0,y) dy \right) \rho(x) dx \leq \]
\[
\leq C_{24} \E \int_{\R^d}\left(\int_{\R^d} K(t,x-y) \f^2(0,y) dy \right) \rho(x) dx \leq C_{24} \|\rho\|_{\infty} \E \|\f(0, \cdot)\|_{B_0^{\rho}}^2<\infty,
\]
where $K$ is the heat kernel in $\R^d$. The estimates for $I_2$ and $I_3$ can be estimated along the lines of Theorem 1 \cite{MisStaYip}.  \\

In order to estimate the second term in (\ref{3.26}), once again we consider two cases: $t \in [0,h]$ and $t \geq h$. If $t \in [0,h]$,
\begin{multline*}
\E \|u_t\|^2_{B_1^\rho} = \E  \int_{-h}^0 \|u(t+\theta)\|_{B_0^\rho}^2 d \theta \leq \E \int_{-h}^0 \|u(s)\|^2_{B_0^\rho} ds + \E \int_0^h \|u(s)\|^2_{B_0^\rho} ds \\
\leq \E  \|\f(t, \cdot)\|^2_{B_1^\rho} + h \sup_{t \geq 0} \E \|u(t)\|^2_{B_0^\rho} < \infty.
\end{multline*}

Finally, if $t \geq h$,
\[
\E \|u_t\|^2_{B_1^\rho} = \E  \int_{-h}^0 \|u(t+\theta)\|^2_{B_0^\rho} d \theta \leq \sup_{t \geq 0} \E \|u(t)\|^2_{B_0^\rho} < \infty,
\]
which completes the proof of the Theorem.
\begin{example}
Let $\bar{f}: \R \to \R$ and $\bar{\sigma}: \R  \to \R$ be Lipschitz functions with Lipschitz constants $L$. Define 
\[
f[\varphi]:= \bar{f}\left(\int_{-h}^0 \varphi(\theta) d \theta\right)
\]
and
\[
\sigma[\varphi]:= \bar{\sigma}\left(\int_{-h}^0 \varphi(\theta) d \theta\right)
\]
Then for all $\varphi_1, \varphi_2 \in B_1^{\rho}$ we have
\[
|f[\varphi_1] - f[\varphi_2]| \leq L \int_{-h}^{0}|\varphi_1(\theta) - \varphi_2(\theta)| d \theta.
\]
Hence
\[
\|f[\varphi_1] - f[\varphi_2]\|_{B_0^\rho}^2 \leq  L^2 \int_{\R^d} \left(\int_{-h}^{0}|\varphi_1(\theta) - \varphi_2(\theta)| d \theta \right)^2 \rho dx \leq L^2 h \|\varphi_1 - \varphi_2\|_{B_1^{\rho}}^2. 
\]
Similarly
\[
\|\sigma[\varphi_1] - \sigma[\varphi_2]\|_{B_0^\rho}^2  \leq L^2 h \|\varphi_1 - \varphi_2\|_{B_1^{\rho}}^2.
\]
Thus $f$ and $\sigma$ are examples of Lipschits maps from $B_1^\rho$ to $B_0^{\rho}$, for which the theorems above apply.
\end{example}

\section{Uniqueness of the invariant measure}

This section is devoted to the proof of Theorem \ref{Th:6}.
\begin{proof}
Let $\mathcal{B}$ be the class of  $\mathcal{F}_t$ measurable $B_0$-valued processes $\xi(t)$, such that
\[
\sup_{t \in \R}\E \|\xi(t)\|_{B_0}^2 < \infty. 
\]
Since 
\[
\sup_{t \in \R} \|\xi(t)\|_B^2 \leq (1+h) \sup_{t \in \R}\E \|\xi(t)\|^2_{B_0}, 
\]
we may follow the procedure from \cite{MisStaYip} and define successive approximations as $u^{(0)} = 0$ and
\begin{equation}\label{iter}
d u^{(n+1)} = (A u^{(n+1)} + f(u^{(n)}_t)) dt + \sigma(u_t^{(n)}) d W(t)
\end{equation}
Then
\[
\sup_{t \in \R}\E \|f(u^{(n)}_t)\|_{B_0}^2 \leq 2\|f(0)\|_{B_0}^2 + 2 L^2 h^2 \sup_{t \in \R} \E \|u^{(n)}(t)\|^2_{B_0} < \infty.
\]
Similarly, 
\[
\sup_{t \in \R}\E \|\sigma(u^{(n)}_t)\|_{B_0}^2 \leq \infty.
\]
Thus by Theorem 5 \cite{MisStaYip}, equation (\ref{iter}) has the unique solution $u^{(n+1)}(t)$ such that
\[
\sup_{t \in \R} \E \|u^{(n+1)}(t)\|^2_{B_0} < \infty
\]
and therefore
\[
\sup_{t \in \R} \E \|u^{(n+1)}(t) \|_B^2 < \infty.
\]
But
\[
\sup_{t \in \R} \E \|u^{(n)}\|_{B_0}^2  \leq (1+h) \sup_{t\in \R} \E \|u^{(n)}(t)\|^2_{B_0} \leq C + h L^2\left(\frac{4}{\lambda_1^2} + \frac{2a}{\lambda_1}\right)\sup_{t \in \R} \E \|u^{(n-1)}\|^2_{B_0}.
\]
Hence for
\begin{equation}\label{smallness}
h L^2\left(\frac{4}{\lambda_1^2} + \frac{2a}{\lambda_1}\right)<1
\end{equation}
in a similar way to \cite{MisStaYip} we can argue that the sequence is in fact Cauchy, and there exists a unique $ u^*(t) $ such that 
\[
\sup_{t \in \R} \E \|u^*(t)\|_{B} < \infty
\]
and
\[
\sup_{t\in\R} \E \|u^{n}(t)-u^*(t)\|^2_B \to 0, n \to \infty.
\]
Furthermore, we can argue that $u^*$ satisfies
\begin{equation}\label{locSol}
u^*(t) = S(t-t_0) u^*(t_0) + \int_{t_0}^t S(t-t_0) f(u^*_s) ds + \int_{t_0}^t S(t-s) \sigma(u^*_s) dW(s).
\end{equation}
Consider any other solution (\ref{locSol}) such that $\eta(t_0)$ $\mathcal{F}_{t_0}$-measurable, and $\E \|\eta(t_0)\|^2_{B}< \infty$. Here $\eta_{t_0} = \varphi(\theta, x)$ is defined on $[-h,0]$. 
Let us show that the solution $\eta$ converges to $u^*$ exponentially. Since we are interested in the behavior of the solutions for large $t$, suppose $t>t_0+h$. Then $t+\theta > t_0$ and $\eta(t)$ is defined via the formula (\ref{locSol}).  Hence
\begin{multline}
\E \|u^*(t) - \eta(t)\|^2_{B_0} \leq 3e^{-\lambda_1(t-t_0)} \E\|u^*(t_0) - \eta(t_0)\|^2_{B_0} + 3 \frac{L^2}{\lambda_1} \int_{t_0}^{t} e^{-\lambda_1(t-s)} \E \|u^*(s) - \eta(s)\|^2_{B_1} ds +  \\
+ 3 L^2 a \int_{t_0}^t e^{-\lambda_1 (t-s)} \E \|u^*(s) - \eta(s)\|^2_{B_1} ds = 3e^{-\lambda_1(t-t_0)} \E\|u^*(t_0) - \eta(t_0)\|^2_{B_0} \\
+ 3 \left(\frac{L^2}{\lambda_1} + L^2 a \right) \int_{t_0}^{t} e^{-\lambda_1(t-s)} \E \|u^*(s) - \eta(s)\|^2_{B_1} ds.
\end{multline}
In addition,
\begin{align*}
&  \E \|u^*(t) - \eta(t) \|^2_{B_1} = \int_{-h}^0 \E \|u^*(t+\theta) - \eta(t+\theta)\|^2_{B_0} d\theta  \leq \\
&  3 \int_{-h}^0 e^{-\lambda_1(t + \theta - s)} \E\|u^*(t_0) - \eta(t_0)\|^2_{B_0} d\theta + 3 \int_{-h}^0 \left(\frac{L^2}{\lambda_1} \int_{t_0}^{t+\theta} e^{-\lambda(t+\theta -s)} \E\|u^*(s) - \eta(s) \|^2_{B_1} ds\right) d \theta + \\
&  +  3 \int_{-h}^0 \left(L^2 a \int_{t_0}^{t+\theta} e^{-\lambda(t+\theta -s)} \E\|u^*(s) - \eta(s) \|^2_{B_1} ds\right) d \theta.
\end{align*}
However, \[
e^{\lambda_1 (t+\theta - s)} \leq e^{-\lambda_1(t-t_0)} \cdot e^{\lambda_1 h},
\]
thus
\begin{multline*}
 \E \|u^*(t) - \eta(t) \|^2_{B_1} \leq 3 h e^{\lambda_1 h}  e^{-\lambda_1(t-t_0)} E \|u^*(t_0) - \eta(t_0)\|^2_{B_0}  \\
 +3 e^{\lambda_1 h} h \left(\frac{L^2}{\lambda_1} + L^2 a \right) \int_{t_0}^t e^{-\lambda_1(t-s)} \E \|u^*(s) - \eta(s)\|^2_{B_1} ds.
\end{multline*}
Altogether,
\begin{multline}
\E \|u^*(t) - \eta(t)\|^2_B \leq (3 e^{\lambda_1 h} h + 3) e^{-\lambda_1(t-t_0)} \E \|u^*(t_0) - \eta(t_0)\|^2_B + \\
\left(3 + 3 h e^{\lambda_1 h}\right) \left(\frac{L^2}{\lambda_1} + L^2 a \right) \int_{t_0}^t e^{-\lambda_1(t-s)} \E \|u^*(s) - \eta(s)\|^2_{B} ds.
\end{multline}
Therefore, if
\begin{equation}\label{smallness}
\left(3 + 3 h e^{\lambda_1 h}\right) \left(\frac{L^2}{\lambda_1} + L^2 a \right):= \gamma_0 L^2 < \lambda_1
\end{equation}
we have
\[
\E \|u^*(t) - \eta(t)\|^2_B \leq (3 e^{\lambda_1 h} h + 3) e^{(\gamma_0 L^2 - \lambda_1)(t-t_0)} \E \|u^*(t_0) - \eta(t_0)\|^2_B.
\]
Then the existence and uniqueness of invariant measure can be established in the same manner as in \cite{MisStaYip}. 
\end{proof}


\begin{thebibliography}{10}

\bibitem{OmaGur}
J.~F.~M. Al-Omari and S.~A. Gourley.
\newblock A nonlocal reaction-diffusion model for a single species with stage
  structure and distributed maturation delay.
\newblock {\em European J. Appl. Math.}, 16(1):37--51, 2005.

\bibitem{AssMan01}
S.~Assing and R.~Manthey.
\newblock Invariant measures for stochastic heat equations with unbounded
  coefficients.
\newblock {\em Stochastic Process. Appl.}, 103(2):237--256, 2003.

\bibitem{ButSch}
O.~Butkovsky and M.~Scheutzow.
\newblock Invariant measures for stochastic functional differential equations.
\newblock {\em Electron. J. Probab.}, 22, 2017.

\bibitem{Car}
Ye. Carkov.
\newblock {\em Random Perturbations of Functional Differential Equations}.
\newblock Riga, 1989.

\bibitem{Cer99}
S.~Cerrai.
\newblock Differentiability of {M}arkov semigroups for stochastic
  reaction-diffusion equations and applications to control.
\newblock {\em Stochastic Process. Appl.}, 83(1):15--37, 1999.

\bibitem{Cho-Mic}
Anna Chojnowska-Michalik.
\newblock Representation theorem for general stochastic delay equations.
\newblock {\em Bull. Acad. Polon. Sci. S\'{e}r. Sci. Math. Astronom. Phys.},
  26(7):635--642, 1978.

\bibitem{Chow}
P.~Chow.
\newblock {\em Stochastic partial differential equations}.
\newblock Chapman \& Hall/CRC Applied Mathematics and Nonlinear Science Series.
  Chapman \& Hall/CRC, Boca Raton, FL, 2007.

\bibitem{CuiYanSun}
J.~Cui, L.~Yan, and X.~Sun.
\newblock Existence and stability for stochastic partial differential equations
  with infinite delay.
\newblock {\em Abstr. Appl. Anal.}, 2014.

\bibitem{DapZab92}
G.~Da~Prato and J.~Zabczyk.
\newblock {\em Stochastic equations in infinite dimensions}, volume~44 of {\em
  Encyclopedia of Mathematics and its Applications}.
\newblock Cambridge University Press, Cambridge, 1992.

\bibitem{DapZab96}
G.~Da~Prato and J.~Zabczyk.
\newblock {\em Ergodicity for infinite-dimensional systems}, volume 229 of {\em
  London Mathematical Society Lecture Note Series}.
\newblock Cambridge University Press, Cambridge, 1996.

\bibitem{HaiMatSch11}
M.~Hairer, J.~Mattingly, and M.~Scheutzow.
\newblock Asymptotic coupling and a general form of {H}arris' theorem with
  applications to stochastic delay equations.
\newblock {\em Probab. Theory Related Fields}, 149(1-2):223--259, 2011.

\bibitem{Hal77}
J.~Hale.
\newblock {\em Theory of functional differential equations}.
\newblock Springer-Verlag, New York-Heidelberg, second edition, 1977.
\newblock Applied Mathematical Sciences, Vol. 3.

\bibitem{IvaKazSwi}
A.~Ivanov, Y.~Kazmerchuk, and A.~Swishchuk.
\newblock Theory, stochastic stability and applications of stochastic delay
  differential equations: a survey of results.
\newblock {\em Differential Equations Dynam. Systems}, 11(1-2):55--115, 2003.

\bibitem{KenStaTsu}
K.~Kenzhebaev, A.~Stanzhytskyi, and A.~Tsukanova.
\newblock Existence and uniqueness results, and the markovian property of
  solutions for a neutral delay stochastic reaction-diffusion equation in
  entire space.
\newblock {\em Dynam. Systems Appl.}, 28(1):19--46, 2019.

\bibitem{KryBog}
N.~Kryloff and N.~Bogoliouboff.
\newblock La th\'eorie g\'en\'erale de la mesure dans son application \`a
  l'\'etude des syst\`emes dynamiques de la m\'ecanique non lin\'eaire.
\newblock {\em Ann. of Math. (2)}, 38(1):65--113, 1937.

\bibitem{LadSolUra}
O.~Ladyzhenskaya, V.~Solonnikov, and N.~Uraltseva.
\newblock {\em Linear and quasilinear equations of parabolic type}, volume~23
  of {\em Translations of mathematical monographs}.
\newblock AMS, Providence, RI, 1968.

\bibitem{ManZau99}
R.~Manthey and T.~Zausinger.
\newblock Stochastic evolution equations in {$L^{2\nu}_\rho$}.
\newblock {\em Stochastics Stochastics Rep.}, 66(1-2):37--85, 1999.

\bibitem{MisStaYip}
O.~Misiats, O.~Stanzhytskyi, and N.~Yip.
\newblock Existence and uniqueness of invariant measures for stochastic
  reaction-diffusion equations in unbounded domains.
\newblock {\em J. Theor. Probab.}, 29(3):996--1026, 2016.

\bibitem{MisStaYip3}
O.~Misiats, O.~Stanzhytskyi, and N.~Yip.
\newblock Asymptotic behavior and homogenization of invariant measures.
\newblock {\em Stochastics and Dynamics}, 19(02), 2019.

\bibitem{MisStaYip2}
O.~Misiats, O.~Stanzhytskyi, and N.~Yip.
\newblock Invariant measures for reaction-diffusion equations with weakly
  dissipative nonlinearities.
\newblock {\em Stochastics}, DOI: 10.1080/17442508.2019.1691212, 2019.

\bibitem{OvsPel}
R.~Ovsepian and A.~Pelczy\'{n}ski.
\newblock On the existence of a fundamental total and bounded biorthogonal
  sequence in every separable {B}anach space, and related constructions of
  uniformly bounded orthonormal systems in {$L^{2}$}.
\newblock {\em Studia Math.}, 54(2):149--159, 1975.

\bibitem{Paz}
A.~Pazy.
\newblock {\em Semigroups of linear operators and applications to partial
  differential equations}, volume~44 of {\em Applied Mathematical Sciences}.
\newblock Springer-Verlag, New York, 1983.

\bibitem{SamMahSta08}
A.~Samoilenko, N.~I. Mahmudov, and A.~N. Stanzhitskii.
\newblock Existence, uniqueness, and controllability results for neutral
  {FSDES} in {H}ilbert spaces.
\newblock {\em Dynam. Systems Appl.}, 17(1):53--70, 2008.

\bibitem{Jason}
S.~Sayed, O.~Misiats, and J.~Clark.
\newblock Electrical control of effective mass, damping, and stiffness of mems
  devices.
\newblock {\em IEEE Sensors Journal}, 17(5):1363--1372, 2017.

\bibitem{StaMogTsu19}
O.~Stanzhytskyi, V.~Mogilova, and A.~Tsukanova.
\newblock On comparison results for neutral stochastic differential equations
  of reaction-diffusion type.
\newblock {\em Modern Mathematics and Mechanics: Understanding Complex
  Systems}, pages 351--395, 2019.

\bibitem{Tan02}
T.~Taniguchi, K.~Liu, and A.~Truman.
\newblock Existence, uniqueness, and asymptotic behavior of mild solutions to
  stochastic functional differential equations in hilbert spaces.
\newblock {\em Journal of Differential Equations}, 181:72 -- 91, 2002.

\bibitem{WanFan}
K.~Wang and M.~Fan.
\newblock {\em Phase Space Theory and Their Application of Functional
  Differential Equations}.
\newblock Science Press, Beijing, 2009.

\end{thebibliography}

\end{document}